\documentclass[11pt]{article}
\usepackage{latexsym}
\usepackage{amsfonts}
\makeatletter

\def\@footnotetext#1{\insert\footins{%
\footnotesize
    \interlinepenalty\interfootnotelinepenalty
    \splittopskip\footnotesep
    \splitmaxdepth \dp\strutbox \floatingpenalty \@MM
    \hsize\columnwidth \@parboxrestore
   \edef\@currentlabel{\csname p@footnote\endcsname\@thefnmark}\@makefntext
    {\rule{\z@}{\footnotesep}\ignorespaces
      #1\strut}}}
\def\abstract{\small\quotation{\hskip-\parindent\sc Abstract.}}

\def\classification{\@ifnextchar [{\@xfootnotenext}%
   {\begingroup\let\protect\noexpand
      \xdef\@thefnmark{}\endgroup
    \@footnotetext}}

\title {}

\begin{document}
\classification {{\it 1991 Mathematics
Subject Classification:} Primary 14E09, 14E25; Secondary
13B10, 13B25.\\ 
$\dagger$) Partially supported by 
  CRCG   Grant 25500/301/01.\\
$\ast$) Partially supported by Hong Kong RGC Grant Project 7126/98P.}  

\begin{center}
{\bf \Large   Peak reduction technique 
 in commutative algebra } 

\bigskip

{\bf   Vladimir Shpilrain}$^{\dagger}$ 

\medskip 

 and 
\medskip

 {\bf Jie-Tai Yu}$^{\ast}$
\smallskip

\end{center} 
\medskip 

\begin{abstract}
\noindent  The ``peak reduction" method is a powerful combinatorial 
technique with applications in many different areas of mathematics 
as well as theoretical computer science. It was introduced 
by Whitehead, a famous topologist and group theorist, who used it
to solve an important algorithmic problem concerning automorphisms of a free
group. Since then, this method was used to solve numerous problems in
group theory,   topology, combinatorics, and probably in some other areas 
as well. 

   In this paper, we  give a survey of  what seems to be the first 
applications  
of the peak reduction technique in commutative algebra and affine algebraic
 geometry. Using this technique, we have contributed toward a classification
 of two-variable polynomials having  classified, up to  an automorphism,  
polynomials of the form $ax^n + by^m + \sum_{im+jn \le
mn} c_{ij} x^i y^j$  (i.e.,  polynomials whose 
Newton polygon is  either a triangle or a line segment). 
 This  has several applications to the study of embeddings of algebraic 
curves in the plane.  In particular,  upon combining our 
method with a well-known theorem of Zaidenberg and Lin, we have shown  
that one can decide ``almost" just by inspection whether or not 
a polynomial  fiber ${\{}p(x,y)=0{\}}$ is an irreducible simply 
connected curve in ${\bf C^2}$. 

 Recently, P.Wightwick  used the idea of peak reduction in 
combination with  splice diagrams technique due to D.Eisenbud and  W.Neumann  
to classify {\it all} two-variable  polynomials over ${\bf C}$ up to  an
automorphism. 

 Another application that we present here, yields  a  decomposition of the
group $Aut(K[x, y])$ in a free product with amalgamation. 

\end{abstract}

\bigskip

\noindent {\bf 1. Introduction }
\bigskip

  In this paper, we  present what seems to be the first applications  
of the  ``peak reduction" method due to Whitehead (see  \cite{White} or 
\cite{LS})  in commutative algebra and affine algebraic geometry. 
In general, this method is used  to find  some kind of canonical 
form of a given object $P$ under the action of a given group (or a semigroup) 
$T$  of transformations. 

 The general idea behind this method is rather simple: one chooses the 
{\it complexity} of  an object $P$ one way or another, and declares a
 canonical form of $P$ an object $P'$ whose complexity is minimal 
among all objects $t(P), ~t \in T$. To actually find a  ``canonical model" 
 $P'$ of a given object $P$, one tries to arrange a sequence of sufficiently 
simple transformations so that the complexity of an object decreases 
{\it at every step}. To prove that such an arrangement is possible, one uses 
``peak reduction"; that means, if in some sequence of simple transformations 
the complexity goes up (or 
remains unchanged) before eventually going down, then there must 
be a pair of {\it subsequent} simple transformations in the sequence (a ``peak")
such 
  that  one of them increases the maximum degree (or leaves it 
unchanged), and then the other one decreases it. Then one tries to prove 
that such a peak can always be reduced. 

 In commutative algebra context, objects are polynomials; their complexity 
is their degree; the group of transformations is the group of polynomial 
automorphisms; simple transformations are elementary and 
linear  automorphisms. (An elementary automorphism is a one that changes 
just one variable). 

  More specific  details on this technique are given in the end of the
Introduction and  in Sections 2 and 3. Here we start by  presenting  results
obtained by using this method. We also point out that, although  this 
technique
was originally created for solving an algorithmic problem,  it can be also 
used to obtain some structural descriptions -- we illustrate this in Section 3. 

\smallskip

 Let $ K[x, y]$ be the polynomial algebra in two variables 
 over a field $K$ of characteristic $0$. The main motivation of our 
paper \cite{ShYu99} was the following problem, which is similar to the 
problem considered by Whitehead in the context of a free group: 

\medskip 

\noindent {\bf Problem 1.} Given two polynomials $p, q \in K[x, y]$, 
find out if there is an automorphism of $K[x, y]$ that takes $p$ to $q$. 
\medskip 

In \cite{ShYu99},  we contributed toward a solution of this problem by
establishing  the following 
\medskip

\noindent {\bf Theorem 1.1.} \cite{ShYu99} Let $p(x,y)= ax^n + by^m + 
\sum_{im+jn \le mn} c_{ij} x^i y^j$, $~a, b, c_{ij} \in K, ~i,j>0$; 
$a, b \ne 0$, 
and  $q(x,y)= Ax^r + By^s + 
\sum_{is+jr \le rs} b_{ij} x^i y^j$, $A, B, b_{ij} \in K, ~i,j>0$; 
 $A, B \ne 0$. Suppose that 
$m$ does not divide $n$, $n$ does not divide $m$, $s$ does not
divide  $r$,$~r$ does not divide  $s$,  and 
max$(m,n) \ne \mbox{max}(r,s)$.  Then  there is no 
 automorphism $\alpha \in Aut(K[x, y])$ that takes $p(x,y)$ to $q(x,y)$. 
\medskip

In some special cases, we can  handle those polynomials with $m$ divisible by
$n$ or vice versa. This is possible, for example, if some fiber of a given 
polynomial admits a one-variable polynomial parametrization $x=u(t); ~y=v(t)$:
\medskip 

\noindent {\bf Proposition 1.2.} \cite{ShYu99}  Suppose the fibers 
${\{}p(x,y)=0{\}}$, ${\{}q(x,y)=0{\}}$ of two polynomials 
$p,q \in \mathbf C[x,y]$, admit  one-variable 
polynomial parametrizations. Then one can effectively find out 
(even without knowing the parametrizations) 
if there is an 
automorphism of $\mathbf C[x,y]$ that takes $p$ to $q$.
\medskip 

 In particular, if some fiber of a given 
polynomial is an irreducible simply 
connected curve, then, by  a well-known theorem of 
Zaidenberg and Lin \cite{ZL}, this  fiber admits a one-variable 
polynomial parametrization. More precisely, they prove 
 that (in case $K= {\mathbf C}$) every polynomial like
that has  a ``canonical model" of the form $x^k - y^l$ with $(k,l)=1$. 
(This generalizes earlier result of  Abhyankar and Moh
\cite{AbMo}). 
Upon combining this with our method, we have the following 
\medskip 

\noindent {\bf Theorem 1.3.} \cite{ShYu99}  Let $p(x,y) \in {\mathbf C}[x,y]$ be 
a polynomial whose fiber ${\{}p(x,y)=0{\}}$ is an irreducible simply 
connected curve. Then some automorphism  of ${\mathbf C}[x,y]$ takes 
$p(x,y)$ to $x^k - y^l$ with $(k,l)=1$, and: 
\smallskip 

\noindent {\bf (a)} max$(k,l) \le \mbox{deg}(p(x,y))$;
\smallskip 

\noindent {\bf (b)} either $k$ or  $l$ divides $\mbox{deg}(p(x,y))$; 
\smallskip 

\noindent {\bf (c)} the  Newton polygon of $p(x,y)$  is either 
a triangle 
or  a line segment, i.e., $p(x,y)$  is of the form 
$ax^n + by^m + {\displaystyle 
\sum_{im+jn \le mn}} c_{ij} x^i y^j$, $i, j > 0, ~a, b \ne 0$. 
 If $m$ does not divide $n$, $n$ does not divide $m$, and $m, n
\ne 0$, then $m=k$  or $l$, 
and $n=l$ or $k$, respectively. Otherwise, either $p(x,y)$  is 
linear, or 
the ``leading" 
 part $ax^n + by^m +  {\displaystyle \sum_{im+jn = mn}} c_{ij} x^i y^j$ ~is a 
proper 
power of some other polynomial.
\medskip 

 Thus, in many situations it is possible to rule out polynomials 
without irreducible simply 
connected fibers just by inspection. In any case, by Proposition 1.2,
there is an  effective procedure for deciding if a given 
polynomial fiber is irreducible and simply connected. 
\medskip 

 Recently, Wightwick \cite{Wightwick} used the idea of peak reduction in 
combination with  splice diagrams technique (see \cite{EN}) to completely 
solve  Problem 1 for two-variable  polynomials over ${\bf C}$. The key ingredient
 of her solution is the following 
\medskip 

\noindent {\bf Theorem 1.4.} \cite{Wightwick}, \cite{NW}
 Let $p\in {\bf C}[x,y]$ be a non-constant polynomial and
  $\phi\in Aut({\bf C}[x,y])$.  Then there is a factorization of $\phi$ in 
a product  $\phi_1\phi_2\dots\phi_n$ of elementary and linear automorphisms, 
such  
that, for   $p_i=\phi_1\phi_2\dots\phi_i(p), ~i=1,\dots,n$, and $p_0=p$,  one
 has:
$$  deg(p_0)\ge\dots\ge deg(p_k)\le\dots\le deg(p_n), ~for
  ~some ~0\le k\le n.$$
Moreover, there is at most one $\phi_i$ for which
$deg(p_{i-1})=deg(p_i)$, and if this happens then $deg(p_i)$ is the
minimal degree $deg(p_k)$.
 \medskip 

 This result provides a procedure for finding a 
``canonical model" for a given polynomial $p$ (i.e., an automorphic image 
of $p$ whose  degree  cannot be reduced by any automorphism). 
Indeed, given a polynomial $p$, we check if there is an elementary 
automorphism
that reduces the degree of $p$. This actually amounts to checking 
automorphisms
of the form $\{x \to x, ~y \to \lambda_1 \cdot  y + \lambda_2 \cdot x^k\}$ 
and 
$\{x \to \lambda_1 \cdot x+ \lambda_2  \cdot y^k, ~y \to y \}$, with  $k \le
deg(p)$. To check if an  elementary automorphism like that (for a particular
$k$) can reduce the degree of $p$, one has  to find out if a specific system 
of polynomial equations for  $\lambda_1, \lambda_2$ has a solution. The latter
 can be done by using  Gr\"obner basis technique (see e.g. \cite{AL}). 
If no elementary automorphism can reduce the degree, then we already 
have a canonical model -- this is precisely the point of Theorem 1.4.

 What is left now to complete the solution of Problem 1, is to be able 
to decide whether or not a canonical model for a given  polynomial $p$ 
can be taken to a canonical model for another given  polynomial $q$ 
by a sequence of elementary and linear automorphisms  of ${\bf C}[x,y])$, 
none of
which changes  the degree of $p$. This, again, ultimately comes down to 
deciding whether or not a specific system  of polynomial equations over 
${\bf C}$ has a solution.

 \medskip 

  We now briefly describe our method that was used 
 to prove Theorem 1.1, Proposition 1.2, and Theorem 1.3, 
leaving the details to the following sections. 

 It is a well-known result of Jung and  van der  Kulk that every automorphism
of $ K[x, y]$ is  a product of elementary and linear automorphisms. 
We want to  get a canonical model for a given polynomial by
  finding  a sequence of
elementary and linear  automorphisms that would reduce the degree {\it at 
every 
step}, until it is further irreducible by any elementary 
automorphism. Then this last polynomial, whose degree is 
irreducible, will be  a canonical model. 

 To  arrange that, we use two principal ideas. First, we mimic   
elementary automorphisms of $ K[x, y]$ by ``elementary 
transformations" of $K[t] \times K[t]$. Second, we use 
Whitehead's idea of ``peak reduction" (see e.g. \cite{LS}) to 
arrange a sequence of elementary transformations of 
$K[t] \times K[t]$ so that the maximum  degree would decrease 
at every step. 
This means the following. If at some 
point of a sequence of ET, the maximum degree goes up (or 
remains unchanged) before eventually going down, then there must 
be a pair of {\it subsequent} ET in our sequence (a ``peak")
such 
  that  one of them increases the maximum degree (or leaves it 
unchanged), and then the other one decreases it. We 
show that such a peak can always be reduced.  
This is described in the next Section 2. 

 While the ``peak reduction" always works for elementary 
transformations of $K[t] \times K[t]$, the first part 
(mimicking elementary automorphisms of $ K[x, y]$ by elementary 
transformations of $K[t] \times K[t]$) is where the 
difficulty is. We managed to do that for polynomials of the 
form given in Theorem 1.1, and also for polynomials 
$p(x,y)$ whose fiber ${\{}p(x,y)=0{\}}$ admits a one-variable 
polynomial parametrization $x=u(t); ~y=v(t)$ (i.e., this 
fiber is a rational curve with one place at infinity). 
The latter was used in proving Proposition 1.2 and Theorem 1.3.
\medskip 

 It would be interesting and important to find applications of 
the  peak reduction method  to the study of the group $Aut(K[x, y,z])$, 
or, at least, of the subgroup of this group generated by tame 
automorphisms (those are products of elementary and linear  automorphisms).
 
A progress in this direction could lead to proving some particular 
 automorphisms of $K[x, y,z]$ to be non-tame, thus resolving a well-known
problem due to Nagata (see e.g. \cite{PMCohn}). In particular, one can ask:

\medskip 

\noindent {\bf Problem 2.} Let $p=p(x,y,z) \in K[x, y,z]$. Is it true 
that if the  degree of $p$ can be decreased by a sequence of elementary and
 linear  automorphisms, then  it can also be decreased by a single elementary
automorphism? 
\medskip

To conclude the Introduction, we mention some earlier  results that were
originally  established by different methods, but can be 
 re-proved in a uniform way by using the peak reduction technique:
\medskip 

\noindent {\bf (1)} \cite{Park} An algorithm for deciding whether or not a
given  matrix  from $GL_2(K[x,y])$ is a product of elementary and diagonal
matrices. 
\medskip 

\noindent {\bf (2)} \cite{ShYu97}  An algorithm for deciding whether or not
 a given polynomial from $K[x,y]$ is a {\it coordinate}, i.e., an automorphic
image  of $x$. This was later generalized in \cite{DrYu} to detecting 
coordinates in $K[z][x,y]$. 
\medskip 

\noindent {\bf (3)} A decomposition of the
group $Aut(K[x, y])$ in a free product with amalgamation. Several 
decompositions of this group have been previously known; see our Section 3 
for details. Also in Section 3, we use the peak reduction method to obtain 
a decomposition which is slightly different from  the  previously known ones.\\

 \noindent {\bf 2. Elementary automorphisms and  peak reduction}
\bigskip

 We give here 
  a somewhat more precise statement of a well-known result of Jung and 
van der  Kulk which can be found in 
\cite[Theorem 6.8.5]{PMCohn}:
\medskip 

\noindent {\bf  Proposition 2.1.}
Every automorphism of $K[x,y]$ is 
a product of linear automorphisms and automorphisms of the form 
$x \to x + f(y); ~y \to y$. More precisely, if $(g_1,
g_2)$ is an automorphism of   $K[x,y]$  such 
  that  $~\mbox{deg}(g_1) \ge \mbox{deg}(g_2)$, say, then either $(g_1, g_2)$ 
is a
linear automorphism, or there exists a unique $~\mu \in 
 K^{\ast}$  and a positive integer $~d~$  such 
  that  $~\mbox{deg}(g_1 - \mu g_2^d) < \mbox{deg}(g_1)$. 
\medskip

 Now we are going to consider the direct product $K[t] \times K[t]$ of two 
copies of the one-variable polynomial algebra over 
 $K$, and introduce the following elementary 
transformations (ET) that can be applied to elements of this algebra: 
\smallskip 

\noindent {\bf (ET1)} $(u, v) \longrightarrow (u+ \mu \cdot 
v^k, v)$ for some $\mu \in K^\ast; ~k \ge 2$. 
\smallskip 

\noindent {\bf (ET2)} $(u, v) \longrightarrow 
(u, v+ \mu \cdot u^k)$. 
\smallskip 

\noindent {\bf (ET3)} a non-degenerate linear transformation 
$(u, v) \longrightarrow (a_1u + a_2v, ~b_1u + b_2v);\\
 a_1, a_2, b_1, b_2 \in K$. 
\smallskip 

  One might notice that some of these transformations  are 
redundant, e.g., (ET1)  is a composition of the other ones. There is 
a reason behind that which will be clear a little later. 

 Our proof of Theorem 1.1 was based on the following 
\medskip

\noindent {\bf  Proposition 2.2.} \cite{ShYu99}  For any pair $(u, v) \in K[t] 
\times K[t]$, there is a (perhaps, empty)  sequence of 
elementary transformations that takes $~(u, v)$ ~to some 
$~(\hat{u}, \hat{v})$   ~such  
that:
\smallskip 

\noindent {\bf (i)} the maximum of the degrees 
of polynomials decreases {\it at every step} in this sequence; 
\smallskip 

\noindent {\bf (ii)} the maximum of the degrees 
in $(\hat{u}, \hat{v})$ is irreducible by {\it any} sequence of 
elementary transformations. 
\smallskip 

 Comment to {\bf (i)}: if it happens so that $u$ and  $v$ have 
 the same leading terms, 
 then, perhaps by somewhat abusing the language, 
we say that the transformation $(u, v) \to (u-v, v)$ reduces 
the maximum of the degrees. 
\medskip 

 We give a proof of Proposition 2.2 here as a sample of our technique.
\medskip 

\noindent {\bf Proof.} We shall use the ``peak reduction" method 
to prove this statement. In this context, this means the following. If at 
some 
point of a sequence of ET, the maximum degree goes up (or 
remains unchanged) before eventually going down, then there must 
be a pair of {\it subsequent} ET in our sequence (a ``peak")
such 
  that  one of them increases the maximum degree (or leaves it 
unchanged), and then the other one decreases it. We are going to 
show that such a peak can always be reduced. In other words, {\it 
if the maximum degree can be decreased by a sequence of ET, then 
it can also be decreased by a single ET}. 
To prove   that, we have 
to consider many different cases, but all of them are quite simple.

 Let $(u, v)$ be a pair of polynomials from $K[t] \times K[t]$ 
with, say, $\mbox{deg}(u) \le \mbox{deg}(v)$,  and let $\alpha_1$ 
and $\alpha_2$ be two  subsequent ET applied to $(u, v)$, 
as described in the previous 
paragraph. Consider several cases: 
\smallskip 

\noindent {\bf (1)} $\alpha_1 : (u, v) \longrightarrow (u+ \mu \cdot 
v^k, v)$ for some $\mu \in K^\ast; ~k \ge 2$. 
 
 This $\alpha_1$ strictly
 increases the maximum degree since $\mbox{deg}(u) \le \mbox{deg}(v)$ by the 
assumption.
 Now we have two possibilities for $\alpha_2$ since 
a linear ET cannot decrease the maximum degree in this situation.  
\smallskip 

\noindent {\bf (a)} $\alpha_2 : (u+ \mu \cdot v^k, v)\longrightarrow (u+ \mu 
\cdot v^k, ~v+\lambda (u+ \mu \cdot v^k)^m)$ for some $\lambda  \in K^\ast; 
~m \ge 2$.  
But this obviously {\it increases} the maximum degree, contrary 
to our assumption. 
\smallskip 

\noindent {\bf (b)} $\alpha_2 : (u+ \mu \cdot v^k, v)\longrightarrow
  (u+ \mu \cdot v^k + \lambda \cdot v^m, v)$.  If this $\alpha_2$ 
 decreases the maximum degree, then we should have  
$\mu \cdot v^k = -\lambda \cdot v^m$, in which case $\alpha_2=\alpha_1^{-1}$, 
and the peak reduction is just 
cancelling out $\alpha_1$ and  $\alpha_2$.
\smallskip 

\noindent {\bf (2)} $\alpha_1 : (u, v) \longrightarrow (u, v+ \mu \cdot u^k)$ 
for some $\mu \in K^\ast; ~k \ge 2$. 

  If this $\alpha_1$ increases the maximum degree, this can only 
happen when $\mbox{deg}(v+ \mu \cdot u^k) = \mbox{deg}(u^k)$, 
in which case we argue exactly as in the case (1). However, 
since $\mbox{deg}(u) \le \mbox{deg}(v)$, it might happen that 
this $\alpha_1$ does not change the maximum degree. Then we 
consider two possibilities for $\alpha_2$:
\smallskip 

\noindent {\bf (a)}  $\alpha_2 : (u, v+ \mu \cdot u^k) \longrightarrow
  (u, v+ \mu \cdot u^k + \lambda \cdot u^m)$.  If this $\alpha_2$ 
 decreases the maximum degree, then we should have $m\ge k$. 
 If  $m=k$, then $\alpha_1 \alpha_2$  is equal to a single ET.
If  $m>k$, then, in order for $\alpha_2$ to decrease the maximum 
degree, we must have $\mbox{deg}(v)$  divisible by $\mbox{deg}(u)$, in which 
case $\alpha_2$ alone would decrease the maximum degree of $(u, v)$, i.e., we 
can get rid of $\alpha_1$. 
\smallskip 

\noindent {\bf (b)} $\alpha_2 : (u, v+ \mu \cdot u^k)\longrightarrow (u+ 
\lambda (v+ \mu \cdot u^k)^m, ~v+ \mu 
\cdot u^k)$. But this $\alpha_2$ can only change the degree 
of the first polynomial in the pair, and  this is not where 
the maximum degree was. 
\smallskip 

\noindent {\bf (3)} $\alpha_1$ is linear, i.e., 
$\alpha_1 : (u, v) \longrightarrow (a_1u + a_2v, b_1u + b_2v); 
 ~a_1, a_2, b_1, b_2 \in K$.  
Again, we have 
two possibilities  for $\alpha_2$.
\smallskip 

\noindent {\bf (a)}  $\alpha_2 : (a_1u + a_2v, b_1u + b_2v) \longrightarrow 
 (a_1u + a_2v, ~b_1u + b_2v + \mu (a_1u + a_2v)^k)$. 
 If $k=1$, then   $\alpha_2$ is linear, and  therefore $\alpha_1 \alpha_2$ 
is a single ET. 
If $k>1$, then $\alpha_2$  might 
decrease  the maximum degree, but this can only happen if 
 $a_2=0$, in which case we could   decrease  the maximum degree 
of $(u, v)$ by a single ET of the type 
(ET2).
\smallskip 

\noindent {\bf (b)}  the case where $\alpha_2$ is of the type 
(ET1), is completely similar.
 \smallskip 

 Thus, in each of the considered cases, if there is a ``peak", 
then  we can reduce the number 
of  ET in the sequence. An obvious inductive argument completes 
the proof of Proposition 2.2. $\Box$ \\

 \noindent {\bf 3.  Decomposing the group of polynomial automorphisms} 
\bigskip

 In this section, we give an application of the ``peak reduction" method 
to produce a new decomposition of the group $Aut(K[x, y])$ in a free
 product with amalgamation. 

 There is a decomposition of the group $Aut(K[x, y])$ in 
    a free product with amalgamation due to Shafarevich \cite{Shaf}; 
  see also \cite[Theorem 6.8.6]{PMCohn}, \cite{Dicks},  \cite{Wright1}
  and references thereto. 
  In \cite{ShYu}, we offered a somewhat more peculiar decomposition.

 To describe our new decomposition and to compare it to previously known ones,
we have to introduce some more notation. 
 \smallskip
 
  We denote:
 \smallskip
 
 \noindent -- by $Af$ the group of affine automorphisms 
 of $K[x, y]$;
 \smallskip
 
 \noindent -- by $UT$ the group of {\it upper triangular} automorphisms;  
those are automorphisms of the  form $x \to a x + p(y); ~y \to 
  by +c$, where  $a,b,c \in K; ~p(y) \in K[y]$. 
 \smallskip
 
 \noindent -- by $LT$  the group of {\it lower triangular} automorphisms; 
those are automorphisms of the  form  $x \to  a x +b; ~y \to cy + p(x)$.
 \smallskip

\noindent -- by $TUT$ the group of  ``twisted" 
upper triangular automorphisms; those 
are automorphisms of the  form $x \to a x + p(y); ~y \to 
 by+cx+d$, where  $a,b,c,d \in K; ~p(y) \in K[y]$. 
\smallskip

\noindent -- by $TLT$  the group of ``twisted" 
 lower triangular automorphisms; those are automorphisms of the  
form $x \to a x +by+c; ~y \to d y + p(x)$.

\smallskip

 Note that  $TUT \cap TLT = Af$. 
\smallskip
 
  There are the following well-known  decompositions in a free product with
amalgamation: 
 \begin{equation}
 \label{1}
 Aut(K[x, y]) = Af \ast_{_{Af \cap UT}} UT = Af \ast_{_{Af \cap LT}} LT. 
 \end{equation}
 
  Our  decomposition has a more  symmetric form:
 \medskip
 
\noindent {\bf Theorem 3.1.} $Aut(K[x, y]) = 
TUT \ast_{_{TUT \cap TLT}}TLT = TUT \ast_{Af}TLT$. 
\medskip

  This cannot be claimed as a brand new result; although it
 probably does not appear anywhere else in exactly this form,
 it can be easily deduced from (1), as well as  (1) 
can be easily deduced from Theorem 1.3. 
 Our proof however is new and  basically self-contained; 
 we only use the aforementioned Jung--van der  Kulk theorem, whereas 
 all known proofs of (1) also use Nagao's theorem \cite{o}. 
 Crucial 
for  the  proof is the following Lemma 3.2; we give a proof of this 
lemma here, as another sample of our method. 
\smallskip 

 We are going to distinguish non-linear 
automorphisms within the groups $TUT$ and  $TLT$. To this end, 
we introduce the following elementary 
transformations (ET)  applied to pairs of polynomials from 
$ K[x, y]$:
\smallskip 

\noindent {\bf (E1)} $(u, v) \longrightarrow (u+ a \cdot 
v^k, v)$ for some $a \in K^\ast; ~k \ge 2$. 
\smallskip 

\noindent {\bf (E2)} $(u, v) \longrightarrow 
(u, v+ a \cdot u^k)$. 
\smallskip 

\noindent {\bf (E3)}  non-degenerate  affine transformations  
$(u, v) \longrightarrow (a_1u + a_2v + c_1, ~b_1u + b_2v + c_2); 
 a_1, a_2, b_1, b_2, c_1, c_2 \in K$. 
\smallskip 

 These elementary transformations generate a group which is 
isomorphic to  $ Aut(K[x, y])$,  by the theorem of Jung and 
van der  Kulk. Now comes 
\medskip

\noindent {\bf  Lemma 3.2.} The group generated by transformations
of the type (E1) and (E2) is a free product of the subgroup  
generated by the transformations (E1) and the one generated by 
the transformations (E2). (Note that both these subgroups are 
abelian). 
\smallskip 

\noindent {\bf Proof} is again based on the ``peak reduction" method. 
 In this context, this means the following. If at some 
point of a sequence of ET (of the type (E1)  or  (E2)), 
the maximum degree goes up (or 
remains unchanged) before eventually going down, then there must 
be a pair of {\it subsequent} ET in our sequence (a ``peak")
such 
  that  one of them increases the maximum degree (or leaves it 
unchanged), and then the other one decreases it. We are going to 
show that such a peak can always be reduced.  

 Let $(u, v)$ be a pair of polynomials from $K[x, y]$ 
with, say, $\mbox{deg}(u) \le \mbox{deg}(v)$,  and let $\gamma_1$ 
and $\gamma_2$ be two  subsequent ET applied to $(u, v)$, 
as described in the previous paragraph. Consider several cases: 
\smallskip 

\noindent {\bf (1)} $\gamma_1 : (u, v) \longrightarrow 
(u+ a \cdot v^k, v)$ for some $a \in K^\ast; ~k \ge 2$. 
 
 This $\gamma_1$ strictly
 increases the maximum degree since $\mbox{deg}(u) \le \mbox{deg}(v)$ by the 
assumption. Now we have two possibilities for $\gamma_2$: 
\smallskip 

\noindent {\bf (a)} $\gamma_2 : (u+ \mu \cdot v^k, v) 
\longrightarrow (u+ \mu \cdot v^k, ~v+ 
b(u+ \mu \cdot v^k)^m)$ for some $b \in K^\ast; ~m \ge 2$.  
But this obviously {\it increases} the maximum degree, contrary 
to our assumption. 
\smallskip 

\noindent {\bf (b)} $\gamma_2 : (u+ \mu \cdot v^k, v) 
\longrightarrow
  (u+ \mu \cdot v^k + b \cdot v^m, v)$.  If this $\gamma_2$ 
 decreases the maximum degree, then we should have  
$a \cdot v^k = -b \cdot v^m$, in which case $\gamma_2=\gamma_1^{-1}$,
 and the peak reduction is just 
cancelling out $\gamma_1$ and  $\gamma_2$.
\smallskip 

\noindent {\bf (2)} $\gamma_1 : (u, v) \longrightarrow 
(u, v+ a \cdot u^k)$ for some $a \in K^\ast; ~k \ge 2$. 

  If this $\gamma_1$ increases the maximum degree, this can only 
happen when $\mbox{deg}(v+ a \cdot u^k) = \mbox{deg}(u^k)$, 
in which case we argue exactly as in the case (1). However, 
since $\mbox{deg}(u) \le \mbox{deg}(v)$, it might happen that 
this $\gamma_1$ does not change the maximum degree. Then we 
consider two possibilities for $\gamma_2$:
\smallskip 

\noindent {\bf (a)}  $\gamma_2 : (u, v+ a \cdot u^k) \longrightarrow
  (u, v+ a \cdot u^k + b \cdot u^m)$.  If this $\gamma_2$ 
 decreases the maximum degree, then we should have $m\ge k$. 
 If  $m=k$, then $\gamma_1 \gamma_2$  is equal to a single ET.
If  $m>k$, then, in order for $\gamma_2$ to decrease the maximum 
degree, we must have $\mbox{deg}(v)$  divisible by 
$\mbox{deg}(u)$, in which case $\gamma_2$ alone would decrease 
the maximum degree of $(u, v)$, i.e., we can switch 
 $\gamma_1$ and $\gamma_2$ (note that these transformations commute
since they are both of the type (E2)). 
\smallskip 

\noindent {\bf (b)} $\gamma_2 : (u, v+ a \cdot u^k) 
\longrightarrow (u+ b (v+ a \cdot u^k)^m, ~v+ a \cdot u^k)$. 
But this $\gamma_2$ can only change the degree 
of the first polynomial in the pair, and  this is not where 
the maximum degree was. 
\smallskip 

 Thus,  given a product of ET of the type (E1)  or  (E2), 
we can, after possibly cancelling pairs of  successive ET of  the 
form $\gamma^{-1} \gamma$ and  switching successive ET of the 
same type, get another product of ET (representing, of course, 
the same automorphism as the given one, call it $\tau$), 
where each factor 
decreases  the maximum degree of a pair of polynomials, starting 
     with the pair $(\tau(x), \tau(y))$, and ending up with the 
pair $(x, y)$. At every step, the choice of ET that can 
decreases  the maximum degree, is unique.

 Therefore,  for any automorphism from the group 
generated by transformations of the type (E1) and (E2), 
there is a unique  alternating product of automorphisms from the 
subgroup  generated by the transformations (E1) and the one generated by 
the transformations (E2). This completes the proof of Lemma 3.2.
$\Box$

\noindent 
 Department of Mathematics, The City  College  of New York, New York, 
NY 10031 
\smallskip

\noindent {\it e-mail address\/}: ~shpil@groups.sci.ccny.cuny.edu \\

\noindent Department of Mathematics, The University of Hong Kong, 
Pokfulam Road, Hong Kong 

\smallskip

\noindent {\it e-mail address\/}: ~yujt@hkusua.hku.hk

\end{document}